\documentclass[journal,twocolumn]{IEEEtran}
\usepackage{amsmath,epsfig}
\usepackage{amssymb}
\usepackage{amstext}
\usepackage{amsmath}
\usepackage{txfonts}
\usepackage{graphicx}
\usepackage{color}
\usepackage{epsfig}

\author{Wei Lu and Namrata Vaswani \\Department of Electrical and Computer Engineering, Iowa State University, Ames, IA \\ \{luwei,namrata\}@iastate.edu}
\title{The Wiener-Khinchin Theorem for Non-wide Sense stationary Random Processes}
\begin{document}

\newcommand{\Dnum}{D_{num}}
\newcommand{\pss}{p^{**,i}}
\newcommand{\fr}{f_{r}^i}

\newcommand{\A}{{\cal A}}
\newcommand{\Z}{{\cal Z}}
\newcommand{\B}{{\cal B}}
\newcommand{\R}{{\cal R}}
\newcommand{\reg}{{\cal G}}
\newcommand{\const}{\mbox{const}}

\newcommand{\trace}{\mbox{trace}}

\newcommand{\hsim}{{\hspace{0.0cm} \sim  \hspace{0.0cm}}}
\newcommand{\he}{{\hspace{0.0cm} =  \hspace{0.0cm}}}

\newcommand{\vect}[2]{\left[\begin{array}{cccccc}
     #1 \\
     #2
   \end{array}
  \right]
  }

\newcommand{\matr}[2]{ \left[\begin{array}{cc}
     #1 \\
     #2
   \end{array}
  \right]
  }
\newcommand{\vc}[2]{\left[\begin{array}{c}
     #1 \\
     #2
   \end{array}
  \right]
  }

\newcommand{\gdot}{\dot{g}}
\newcommand{\Cdot}{\dot{C}}
\newcommand{\re}{\mathbb{R}}
\newcommand{\n}{{\cal N}}  
\newcommand{\N}{{\overrightarrow{\bf N}}}  
\newcommand{\chat}{\tilde{C}_t}
\newcommand{\chati}{\chat^i}

\newcommand{\cmin}{C^*_{min}}
\newcommand{\twi}{\tilde{w}_t^{(i)}}
\newcommand{\twj}{\tilde{w}_t^{(j)}}
\newcommand{\wi}{{w}_t^{(i)}}
\newcommand{\twio}{\tilde{w}_{t-1}^{(i)}}

\newcommand{\tWi}{\tilde{W}_n^{(m)}}
\newcommand{\tWj}{\tilde{W}_n^{(k)}}
\newcommand{\Wi}{{W}_n^{(m)}}
\newcommand{\tWio}{\tilde{W}_{n-1}^{(m)}}

\newcommand{\ds}{\displaystyle}

\newcommand{\SAR}{S$\!$A$\!$R }
\newcommand{\MAR}{MAR}
\newcommand{\MMRF}{MMRF}
\newcommand{\AR}{A$\!$R }
\newcommand{\GMRF}{G$\!$M$\!$R$\!$F }
\newcommand{\DTM}{D$\!$T$\!$M }
\newcommand{\MSE}{M$\!$S$\!$E }
\newcommand{\RCS}{R$\!$C$\!$S }
\newcommand{\uomega}{\underline{\omega}}
\newcommand{\y}{v}
\newcommand{\x}{w}
\newcommand{\lu}{\mu}
\newcommand{\g}{g}
\newcommand{\s}{{\bf s}}
\newcommand{\bft}{{\bf t}}
\newcommand{\refmap}{{\cal R}}
\newcommand{\totrefl}{{\cal E}}
\newcommand{\beq}{\begin{equation}}
\newcommand{\eeq}{\end{equation}}
\newcommand{\bdm}{\begin{displaymath}}
\newcommand{\edm}{\end{displaymath}}
\newcommand{\hatz}{\hat{z}}
\newcommand{\hatu}{\hat{u}}
\newcommand{\tilz}{\tilde{z}}
\newcommand{\tilu}{\tilde{u}}
\newcommand{\hhatz}{\hat{\hat{z}}}
\newcommand{\hhatu}{\hat{\hat{u}}}
\newcommand{\tilc}{\tilde{C}}
\newcommand{\hatc}{\hat{C}}
\newcommand{\tim}{n}

\newcommand{\ssp}{\renewcommand{\baselinestretch}{1.0}}
\newcommand{\defd}{\mbox{$\stackrel{\mbox{$\triangle$}}{=}$}}
\newcommand{\goes}{\rightarrow}
\newcommand{\tends}{\rightarrow}
\newcommand{\defn}{\triangleq} 
\newcommand{\se}{&=&}
\newcommand{\sdefn}{& \defn  &}
\newcommand{\sle}{& \le &}
\newcommand{\sge}{& \ge &}
\newcommand{\plusminus}{\stackrel{+}{-}}
\newcommand{\Ey}{E_{Y_{1:t}}}
\newcommand{\ey}{E_{Y_{1:t}}}

\newcommand{\equivto}{\mbox{~~~which is equivalent to~~~}}
\newcommand{\nonzero}{i:\pi^n(x^{(i)})>0}
\newcommand{\nonzeroc}{i:c(x^{(i)})>0}

\newcommand{\supn}{\sup_{\phi:||\phi||_\infty \le 1}}
\newtheorem{theorem}{Theorem}
\newtheorem{lemma}{Lemma}
\newtheorem{corollary}{Corollary}
\newtheorem{definition}{Definition}
\newtheorem{remark}{Remark}
\newtheorem{example}{Example}
\newtheorem{ass}{Assumption}
\newtheorem{fact}{Fact}
\newtheorem{heuristic}{Heuristic}
\newcommand{\eps}{\epsilon}
\newcommand{\bd}{\begin{definition}}
\newcommand{\ed}{\end{definition}}
\newcommand{\udq}{\underline{D_Q}}
\newcommand{\td}{\tilde{D}}
\newcommand{\epsinv}{\epsilon_{inv}}
\newcommand{\al}{\mathcal{A}}

\newcommand{\bfx} {\bf X}
\newcommand{\bfy} {\bf Y}
\newcommand{\bfz} {\bf Z}
\newcommand{\ddas}{\mbox{${d_1}^2({\bf X})$}}
\newcommand{\ddbs}{\mbox{${d_2}^2({\bfx})$}}
\newcommand{\dda}{\mbox{$d_1(\bfx)$}}
\newcommand{\ddb}{\mbox{$d_2(\bfx)$}}
\newcommand{\xinc}{{\bfx} \in \mbox{$C_1$}}
\newcommand{\eqa}{\stackrel{(a)}{=}}
\newcommand{\eqb}{\stackrel{(b)}{=}}
\newcommand{\eqe}{\stackrel{(e)}{=}}
\newcommand{\leqc}{\stackrel{(c)}{\le}}
\newcommand{\leqd}{\stackrel{(d)}{\le}}

\newcommand{\leqa}{\stackrel{(a)}{\le}}
\newcommand{\leqb}{\stackrel{(b)}{\le}}
\newcommand{\leqe}{\stackrel{(e)}{\le}}
\newcommand{\leqf}{\stackrel{(f)}{\le}}
\newcommand{\leqg}{\stackrel{(g)}{\le}}
\newcommand{\leqh}{\stackrel{(h)}{\le}}
\newcommand{\leqi}{\stackrel{(i)}{\le}}
\newcommand{\leqj}{\stackrel{(j)}{\le}}

\newcommand{\w}{{W^{LDA}}}
\newcommand{\halpha}{\hat{\alpha}}
\newcommand{\hsigma}{\hat{\sigma}}
\newcommand{\slmax}{\sqrt{\lambda_{max}}}
\newcommand{\slmin}{\sqrt{\lambda_{min}}}
\newcommand{\lmax}{\lambda_{max}}
\newcommand{\lmin}{\lambda_{min}}

\newcommand{\da} {\frac{\alpha}{\sigma}}
\newcommand{\chka} {\frac{\check{\alpha}}{\check{\sigma}}}
\newcommand{\sumo}{\sum _{\underline{\omega} \in \Omega}}
\newcommand{\distance}{d\{(\hatz _x, \hatz _y),(\tilz _x, \tilz _y)\}}
\newcommand{\col}{{\rm col}}
\newcommand{\rcs}{\sigma_0}
\newcommand{\CalR}{{\cal R}}
\newcommand{\df}{{\delta p}}
\newcommand{\dq}{{\delta q}}
\newcommand{\dZ}{{\delta Z}}
\newcommand{\pprime}{{\prime\prime}}

\newcommand{\vn}{N}

\newcommand{\bv}{\begin{vugraph}}
\newcommand{\ev}{\end{vugraph}}
\newcommand{\bi}{\begin{itemize}}
\newcommand{\ei}{\end{itemize}}
\newcommand{\ben}{\begin{enumerate}}
\newcommand{\een}{\end{enumerate}}
\newcommand{\be}{\protect\[}
\newcommand{\ee}{\protect\]}
\newcommand{\bean}{\begin{eqnarray*} }
\newcommand{\eean}{\end{eqnarray*} }
\newcommand{\bea}{\begin{eqnarray} }
\newcommand{\eea}{\end{eqnarray} }
\newcommand{\nn}{\nonumber}
\newcommand{\ba}{\begin{array} }
\newcommand{\ea}{\end{array} }
\newcommand{\ep}{\mbox{\boldmath $\epsilon$}}
\newcommand{\epp}{\mbox{\boldmath $\epsilon '$}}
\newcommand{\Lep}{\mbox{\LARGE $\epsilon_2$}}
\newcommand{\und}{\underline}
\newcommand{\pdif}[2]{\frac{\partial #1}{\partial #2}}
\newcommand{\odif}[2]{\frac{d #1}{d #2}}
\newcommand{\dt}[1]{\pdif{#1}{t}}
\newcommand{\urho}{\underline{\rho}}

\newcommand{\spc}{{\cal S}}
\newcommand{\tspc}{{\cal TS}}

\newcommand{\uv}{\underline{v}}
\newcommand{\us}{\underline{s}}
\newcommand{\uc}{\underline{c}}
\newcommand{\utheta}{\underline{\theta}^*}
\newcommand{\ualpha}{\underline{\alpha^*}}

\newcommand{\uxy}{\underline{x}^*}
\newcommand{\uxyj}{[x^{*}_j,y^{*}_j]}
\newcommand{\arcl}[1]{arclen(#1)}
\newcommand{\one}{{\mathbf{1}}}

\newcommand{\uxyjt}{\uxy_{j,t}}
\newcommand{\E}{\mathbb{E}}

\newcommand{\rhomat}{\left[\begin{array}{c}
                        \rho_3 \ \rho_4 \\
                        \rho_5 \ \rho_6
                        \end{array}
                   \right]}
\newcommand{\deltat}{\tau} 
\newcommand{\deltatt}{\Delta t_1}
\newcommand{\ceil}[1]{\ulcorner #1 \urcorner}

\newcommand{\xxi}{x^{(i)}}
\newcommand{\txi}{\tilde{x}^{(i)}}
\newcommand{\txj}{\tilde{x}^{(j)}}

\newcommand{\mi}[1]{{#1}^{(m,i)}}

\date{}
\maketitle
\begin{abstract}
We extend the Wiener-Khinchin theorem to non-wide sense stationary (WSS) random processes, i.e. we prove that, under certain assumptions, the power spectral density (PSD) of any random process is equal to the Fourier transform of the time-averaged autocorrelation function. We use the theorem to show that bandlimitedness of the PSD implies bandlimitedness of the generalized-PSD for a certain class of non-WSS signals. This fact allows us to apply the Nyquist criterion derived by Gardner for the generalized-PSD.
\end{abstract}

\section{Introduction}
The Power Spectral Density (PSD) defined in (\ref{defpsd2}) of a random process is the expected value of its normalized periodogram, with the duration over which the periodogram is computed approaching infinity. For wide sense
stationary (WSS) processes, the Wiener-Khinchin theorem
\cite{YatesGoodman} shows that the PSD
is equal to the Fourier transform of the autocorrelation function(treated as a function of the delay). 
 Wiener-Khinchin theorem is a very fundamental result because it can be used to perform spectral analysis of WSS random processes whose Fourier transform may not exist. 

In this work we give a complete rigorous proof of the nonstationary
analog of the Wiener-Khinchin theorem, i.e. we show that
under certain assumptions, the PSD of a random process, defined in (1), is
equal to the Fourier transform of the time-averaged autocorrelation
function.

While similar ideas have been introduced in some textbooks on random
processes\cite{CooperMcGillem}\cite{Gardnerbook}\cite{commsyseng}\cite{pebble}, (with the aim of generalizing the stationary case), with the exception of \cite{pebble}, the results are mostly incomplete. We discuss this in section 2.1 and give details in \cite{webref}. Also, the result of \cite{pebble} uses a different
sufficient condition than ours. Note
that the work of \cite{Generalization of the Wiener-Khinchin
theorem} has a very related title, but a completely different
contribution from the current work. It proves Wiener-Khinchin theorem for a generalization of autocorrelation.

There are other commonly used definitions of the PSD for non-WSS
processes, e.g. the generalized-PSD \cite{samplenonstat},
$K(u,v)$, which is the 2D Fourier transform of the autocorrelation or the
evolutionary spectral density function \cite{PriestleyRao}. In
\cite{spectralspatial}, $K(u,v)$ has been shown to be equal to the
covariance between the Fourier transform coefficients of the signal at $u$ and at
$v$.

Our result is important because one can now do 1-dimensional
spectral analysis of certain types of non-WSS signals using the PSD.
The generalized-PSD,  $K(u,v)$, defines a 2D Fourier transform which is much more
expensive to estimate (as explained in section 3) and more difficult to
interpret. 
We prove in section 3 that at least for certain types of non-WSS signals, the
bandwidth computed using the PSD can be used for subsampling the
signal. This is done by using our result to show that
bandlimitedness of the PSD implies bandlimitedness of the
generalized-PSD for a certain class of non-WSS signals and then
using Gardner's result \cite{samplenonstat} for subsampling.
One motivating application for the above is analyzing
piecewise stationary signals for which the boundaries between pieces
are not known. It is computationally more efficient if one can first
uniformly subsample the signal to perform dimension reduction (by
using the bandwidth computed from the estimated PSD) and then use
existing techniques to find the piece boundaries or to perform
inference tasks such as signal classification.


A common practical approach for computing the PSD of a WSS process
is to break up a single long sequence into pieces and to use each
piece as a different realization of the process. But this cannot be
done for a non-WSS process. Multiple realizations are required
before either the autocorrelation or the PSD can be estimated and
these are often difficult to obtain in practice. One application where multiple realizations are available is
to analyze time sequences of spatially non-WSS signals, which
are temporally independent and identically distributed (i.i.d) or
temporally stationary and ergodic. For example, if the sequence is a time sequence of
contour deformations or a time sequence of images (2D spatial
signals) when temporal i.i.d-ness or stationarity is a valid
assumption.
%
The contour deformation ``signal" at a given time is a 1D function of
contour arclength and is often spatially piecewise stationary, e.g. very often one region of the contour deforms much more than the
others. Temporal stationarity is a valid assumption in many
practical applications such as when analyzing human body contours
for gait recognition or analyzing brain tumor contour deformations
and thus the time sequence can be used to compute the expectations.

{\em Paper Organization: } In Sec. 2, we give the complete proof of
Wiener-Khinchin for non-WSS processes. In Sec. 3, sampling theory for
a particular class of non-WSS processes is discussed. Application to
time sequences of spatially
non-WSS signals is shown in Sec. 4. Conclusions are given in Sec. 5.%

\setlength{\arraycolsep}{0.04cm}

\section{Non-WSS Wiener-Khinchin Theorem}
The PSD of any random process, $x(t)$, is defined as
\begin{equation}
S_x(\omega) \triangleq \lim_{T\rightarrow
\infty}\frac{E[|X_T(\omega)|^2]}{2T},\ \ X_T(\omega) \defn
\int_{-T}^Tx(t)e^{-j\omega t}dt \label{defpsd2}
\end{equation}
$E[.]$ denotes expectation w.r.t. the pdf of the process $x(t)$.
Let $R_{x}(t_1,t_2) \defn E[x(t_1)x(t_2)]$ be the autocorrelation of
$x(t)$.
Assume%
\bea \int_{-T}^T\int_{-T}^T E[|x(t_1) x(t_2)|] dt_1 dt_2 <\infty
\label{sigmamax} \eea for any finite $T$\footnote{A sufficient condition for (2) is: there exists a $C < \infty$ such that $ R_{|x|}(t_1,t_2) < C$ almost everywhere (except on a set of measure zero)}.
Because of (2), Fubini's theorem \cite[Chap 12]{royden} can be applied to move the expectation inside the integral and to change the order of the integrals in (1). Also, assume that the ``maximum absolute autocorrelation
function", \bea
R_x^m(\tau)\triangleq \sup_{t}|R_{x}(t,t-\tau)| 
\label{rxmax2}
\end{eqnarray}
is integrable.
%
%
Moving the expectation inside the integrals in (\ref{defpsd2}) and defining $\tau \defn t_2-t_1$, we get
\bea
S_x(\omega) \se \lim_{T\rightarrow \infty}\frac{1}{2T}E[\int_{-T}^Tx(t_1)e^{j\omega t_1}dt_1\int_{-T}^Tx(t_2)e^{-j\omega t_2}dt_2] \nn \\
\se \lim_{T\rightarrow
\infty}\frac{1}{2T}\int_{-T}^T\int_{t_2-T}^{t_2+T}R_{x}(t_2,t_2-\tau)e^{-j\omega\tau}d\tau
\ dt_2 \eea Changing the integral order and the corresponding limits
of integration, in a fashion similar to the proof of the stationary
case \cite{YatesGoodman},
\begin{eqnarray}
S_x(\omega)=\lim_{T\rightarrow \infty} \frac{1}{2T}\int_{\tau=0}^{2T}\int_{t_2=\tau-T}^TR_{x}(t_2,t_2-\tau)e^{-j\omega \tau}dt_2d\tau+ \nonumber\\
\lim_{T\rightarrow \infty}\frac{1}{2T}\int_{\tau=-2T}^0\int^{\tau+T}_{t_2=-T}R_{x}(t_2,t_2-\tau)e^{-j\omega \tau}dt_2d\tau \nonumber  \\
=\lim_{T\rightarrow \infty} \frac{1}{2T}\int_{\tau=-2T}^{2T}\int_{t_2=-T}^TR_{x}(t_2,t_2-\tau)e^{-j\omega \tau}dt_2d\tau- \nonumber\\
\lim_{T\rightarrow \infty} \frac{1}{2T}\int_{\tau=0}^{2T}\int_{t_2=-T}^{\tau-T}R_{x}(t_2,t_2-\tau)e^{-j\omega \tau}dt_2d\tau- \nonumber\\
\lim_{T\rightarrow \infty}\frac{1}{2T}\int_{\tau=-2T}^0\int^{T}_{t_2=\tau+T}R_{x}(t_2,t_2-\tau)e^{-j\omega \tau}dt_2d\tau \label{defsx2}
\end{eqnarray}
\begin{displaymath}
\text{Define~~~~~~~~} g_{T}(\tau)\triangleq \left\{\begin{array}{ll}
1 & \textrm{if $-2T \le \tau \le 2T$}\\
0 & \textrm{otherwise}
\end{array} \right.
\end{displaymath}
Then the first term of (\ref{defsx2}) becomes
\bea   
\lim_{T\rightarrow \infty} \frac{1}{2T}\int_{\tau=-2T}^{2T}\int_{t_2=-T}^TR_{x}(t_2,t_2-\tau)e^{-j\omega \tau}dt_2d\tau \nonumber \\= \lim_{T\rightarrow \infty}\int_{-\infty}^{\infty} f_T(\tau)d\tau \quad \quad \quad \quad \quad \quad \quad \quad \quad \quad \quad \ \ \eea \\
where $f_{T}(\tau)\triangleq \frac{1}{2T}\int_{-T}^{T}R_x(t_2,t_2-\tau)dt_2 \ g_{T}(\tau) \ e^{-j\omega \tau} $.
It is easy to see \bea |f_{T}(\tau)| \le \frac{2T
}{2T} R_x^m(\tau)g_{T}(\tau)= R_x^m(\tau)g_{T}(\tau)\le R_x^m(\tau)
\eea
Since $R^m_x(\tau)$ is integrable, we can use dominated convergence theorem \cite[Chap 4]{royden} to move the limit inside the outer integral in the first term of (\ref{defsx2}). Thus, this term becomes
\begin{eqnarray}
&\ &\int_{-\infty}^{\infty}\lim_{T\rightarrow \infty} [\frac{1}{2T}\int_{-T}^{T}R_x(t_2,t_2-\tau)dt_2 g_{T}(\tau)] e^{-j\omega \tau}d\tau \nn \\
&=& \int_{-\infty}^{\infty}\lim_{T\rightarrow \infty}
[\frac{1}{2T}\int_{-T}^{T}R_x(t_2,t_2-\tau)dt_2
] e^{-j\omega  \tau}d\tau \label{defsx3}
\end{eqnarray}
Now, for the second term in (\ref{defsx2}), we know
\begin{eqnarray}
\lim_{T\rightarrow \infty} |\frac{1}{2T}\int_{\tau=0}^{2T}\int_{t_2=-T}^{\tau-T}R_{x}(t_2,t_2-\tau)e^{-j\omega \tau}dt_2d\tau|  \nonumber \\
\le \lim_{T\rightarrow \infty}  \frac{1}{2T}\int_{\tau=0}^{2T}\int_{t=-T}^{\tau-T}|R_{x}(t_2,t_2-\tau)e^{-j\omega \tau}|dt_2 d\tau
\end{eqnarray}
From (\ref{rxmax2}), we know $|R_x(t_2,t_2-\tau)|\le R_x^m(\tau)$. Hence,
\begin{equation}
\lim_{T\rightarrow \infty} |\frac{1}{2T}\int_{\tau=0}^{2T}\int_{t_2=-T}^{\tau-T}R_{x}(t_2,t_2-\tau)e^{-j\omega \tau}dt_2d\tau| \le \lim_{T\rightarrow \infty}  \int_{\tau=0}^{2T} \frac{\tau}{2T} R_{x}^m(\tau)d\tau \label{secondpart1}
\end{equation}
\begin{displaymath}
\text{Now, define~~~~~~~~} \phi_{T}(\tau)\triangleq \left\{\begin{array}{ll}
\frac{\tau}{2T} & \textrm{if $0 \le \tau \le 2T$}\\
0 & \textrm{otherwise}
\end{array} \right.
\end{displaymath}
Then $\lim_{T\rightarrow \infty}  \int_{\tau=0}^{2T} \frac{\tau}{2T} R_{x}^m(\tau)d\tau=\lim_{T\rightarrow \infty}  \int_{\tau=0}^{\infty} \phi_T(\tau) R_{x}^m(\tau)d\tau$. We can apply dominated convergence theorem again since
$|\phi_T(\tau) R_{x}^m(\tau)| \le R_{x}^m(\tau)$ which is integrable. Thus,
this implies
\begin{equation}
\lim_{T\rightarrow \infty}  \int_{\tau=0}^{\infty} \phi_T(\tau) R_{x}^m(\tau)d\tau= \int_{\tau=0}^{\infty} \lim_{T\rightarrow \infty} [\phi_T(\tau)R_{x}^m(\tau)]d\tau=0
\end{equation}
Therefore, (\ref{secondpart1}) forces
\begin{eqnarray}
& & 0 \le \lim_{T\rightarrow \infty} |\frac{1}{2T}\int_{\tau=0}^{2T}\int_{t_2=-T}^{\tau-T}R_{x}(t_2,t_2-\tau)e^{-j\omega \tau}dt_2d\tau| \le 0 \nn \\
& \implies & \lim_{T\rightarrow \infty} \frac{1}{2T}\int_{\tau=0}^{2T}\int_{t_2=-T}^{\tau-T}R_{x}(t_2,t_2-\tau)e^{-j\omega \tau}dt_2d\tau = 0 \label{limzero}
\end{eqnarray}
Thus, the second term of (\ref{defsx2}) is 0. In an analogous fashion, we can show that the third term is also 0. Thus, we finally get \bea \label{nswk}
S_x(\omega) \se \int_{-\infty}^{\infty} \bar{R_x} (\tau) e^{-j\omega \tau} d\tau, \ \
where    \\
 \bar{R_x} (\tau) \sdefn  \lim_{T \tends \infty} \frac{1}{2T}\int_{-T}^{T}R_x(t_2,t_2-\tau)dt_2
 \label{rxbar2}
\eea This is the Wiener-Khinchin result for any general random process, i.e.
\begin{theorem}[Wiener-Khinchin Theorem for Non-WSS Processes]
If (\ref{sigmamax}) holds, and if $R_x^m(\tau)$ defined in (\ref{rxmax2}) is integrable,
then (\ref{nswk}) holds, i.e. the PSD, defined in (\ref{defpsd2}), is
equal to Fourier transform of the ``averaged autocorrelation function", defined in
(\ref{rxbar2}). \label{wk}
\end{theorem}
\subsection{Discussion of related results}
Related ideas are introduced in several textbooks\cite{CooperMcGillem} \cite{Gardnerbook} \cite{commsyseng}\cite{pebble} , but except for \cite{pebble}, the results in the rest of them are either different from ours or incomplete.\\
\cite{CooperMcGillem} does not justify why the limit can be moved inside the integral in equation 7.38.
\cite{Gardnerbook} has a similar problem and also his result says that $S(f)$ is equal to the time average of the instantaneous PSD. The condition of \cite{commsyseng} given on page 179 is not sufficient either. Fourier transform of the averaged autocorrelation function may not exist even if his condition is satisfied. Peeble's result\cite{pebble} gives a different sufficient condition from our result (He assumes absolute integrability of PSD). We discuss the results of the above books in detail in \cite{webref}.

\section{Subsampling A Class of Non-WSS Processes}
We use Gardner's result \cite{samplenonstat} for
subsampling nonstationary random signals along with the PSD, instead
of the generalized-PSD, since the PSD is much less expensive to
compute. For a $T$-length signal for which $P$ realizations are
available, the PSD defined in (1) can be estimated in $O(P T\log T)$ time (need to
estimate $P$ $T$-length Fourier transform and average their square magnitudes).
Generalized-PSD computes $R_x(t_1,t_2)$ using $P$ realizations
(takes $O(PT^2)$ time), followed by computing its 2D Fourier transform (takes
$O(T^2 \log T)$ time), i.e. it requires $O(PT^2 + T^2 \log T)$ time. Clearly $O(PT log T )  < O (PT^2)$.

Gardner's result says that if the generalized PSD of $x$,
$K_x(u,v)$, is bandlimited in both dimensions, $x$ can be
reconstructed exactly (in the mean square sense), from its
uniformly-spaced samples taken at a rate that is higher than twice
the maximum bandwidth in either dimension.
To use this result with the PSD, we need to show that bandlimitedness of the PSD implies bandlimitedness of the generalized-PSD. The most general case for which this can be done will be studied in future. We show it for the following class of non-WSS random signals, which can be used to model many commonly occurring random processes including many piecewise stationary ones. This is one of the four classes of nonstationary processes described in \cite{Nonstatinterp}.%
\bd[Class NS1] A random signal, $x(t)$, belongs to the class NS1 if
it can be represented as the output of  ``nonstationary white
Gaussian noise", $w(t)$, passed through a stable linear time
invariant (LTI) system, denoted $h(t)$, i.e. $x(t) = w(t) \star
h(t)$, where $\star$ denotes convolution; $w(t)$ is a Gaussian
process with $E[w(t)]=0, \ \forall t$, $R_{w}(t,t') = \sigma_w^2(t)
\delta(t-t')$ and $\sigma_{w,max}^2 \defn \max_t \sigma_w^2(t) <
\infty$; and $h(t)$ satisfies $\int_{-\infty}^{\infty} |h(t)|dt <
\infty$. \ed We show below that for this class of signals,
$\bar{R_x}(\tau) = \bar{\sigma_w^2} h(\tau) \star h(-\tau)$ and
that Theorem \ref{wk} can be applied to show that $S_x(\omega) =  \bar{\sigma_w^2}
|H(\omega)|^2$ where $H(\omega)$ denotes the Fourier transform of $h(t)$.\\
Using the definition of NS1 signals, $\bar{R_x}(\tau)$ is written as
\begin{eqnarray}
\bar{R_x}(\tau) \se \lim_{T\rightarrow \infty} \frac{1}{2T}\int_{-T}^{T}\int_{-\infty}^{\infty}\int_{-\infty}^{\infty}h(\tau_1)h(\tau_2)E[w(t_1-\tau_1)  \nn \\
&& w(t_1-\tau-\tau_2)]d\tau_2d\tau_1 dt_1 \nn \\
\se \lim_{T\rightarrow \infty} \frac{1}{2T}
\int_{-T}^{T}\int_{-\infty}^{\infty}h(\tau_1)h(\tau_1-\tau)\sigma_w^2(t_1-\tau_1)d\tau_1dt_1
\ \ \ \ \
\end{eqnarray}
By changing the order of integration and defining $t=t_1-\tau_1$, we
get
\begin{eqnarray}
\bar{R_x}(\tau)=\lim_{T\rightarrow \infty}\int_{-\infty}^{\infty}
\frac{1}{2T}\int_{-T-\tau_1}^{T-\tau_1}\sigma_w^2(t)dt \
h(\tau_1)h(\tau_1-\tau)d\tau_1
\end{eqnarray}
The integrand, $\psi(\tau,\tau_1) \defn
\frac{1}{2T}\int_{-T-\tau_1}^{T-\tau_1}\sigma_w^2(t)dt \
h(\tau_1)h(\tau_1-\tau) \le \sigma^2_{w,max}
h(\tau_1-\tau)h(\tau_1)$. Since $h(\tau_1-\tau)h(\tau_1)$ is
absolutely integrable w.r.t. $\tau_1$ (follows from Cauchy-Schwartz
inequality and the fact that a stable $h(t)$ implies that
$\int_{-\infty}^{\infty} h^2(t) dt < \infty$), we can apply the
dominated convergence theorem \cite{royden} to move the limit inside
to get
\begin{eqnarray}
\bar{R_x}(\tau)=\int_{-\infty}^{\infty} \lim_{T\rightarrow \infty}
\{\frac{1}{2T}\int_{-T-\tau_1}^{T-\tau_1}\sigma_w^2(t)dt\} \
h(\tau_1)h(\tau_1-\tau)d\tau_1
\end{eqnarray}
By using an argument similar to that used in (\ref{limzero}) and the
paragraph below it, we can replace $\lim_{T\rightarrow \infty}
\frac{1}{2T}\int_{-T-\tau_1}^{T-\tau_1} [.]$ by $\lim_{T\rightarrow
\infty} \frac{1}{2T}\int_{-T}^{T} [.]$. Thus, we get \bea
\label{rxbar_ns1}
\bar{R_x}(\tau) \se \int_{-\infty}^{\infty} \bar{\sigma_{w}^2} h(\tau_1)h(\tau_1-\tau)d\tau_1 =  \bar{\sigma_{w}^2} h(\tau) \star h(-\tau), \ \ \ \ \\
where \ \ \bar{\sigma_{w}^2} \sdefn  \lim_{T\rightarrow
\infty}\frac{1}{2T}\int_{-T}^{T}\sigma_w^2(t)dt \eea By taking the
Fourier transform of both sides of (\ref{rxbar_ns1}), and using Theorem \ref{wk}
for the left hand side, we get that \bea S_x(\omega) =
\bar{\sigma_w^2} |H(\omega)|^2 \label{SxNS1} \eea
Theorem \ref{wk}
can be applied because
(a) $R_{x}^m(\tau) =  \sup_{t}
|R_x(t,t-\tau)| \le \sigma_{w,max}^2 |h(\tau) \star h(-\tau)|$ which is
absolutely integrable (since $h(\tau) \star h(-\tau)$ is stable)
and (b) the sufficient condition given in footnote 1 for (\ref{sigmamax}) to hold is satisfied since $R_{|x|}(t_1,t_2)\le \sup_t R_{|x|}(t,t)=\sup_t R_x(t,t) \le \sigma_{w,max}^2 \int_{-\infty}^{\infty}
h^2(\tau_1) d\tau_1\triangleq C$.

From (\ref{SxNS1}), for signals belonging to NS1, if $S_x(\omega)$ is bandlimited, it implies that $H(\omega)$ is bandlimited. This in turn implies that $K_x(u,v) = \int \int R_x(t_1,t_2)e^{-j(ut_1- vt_2)}dt_1 dt_2 = H(u)H^*(v) K_w(u,v)$ is also bandlimited with the same bandwidth. Since $K_x(u,v)$ is bandlimited, Gardner's result \cite{samplenonstat} applies. Thus we have the following corollary.
\begin{corollary}
For random signals $x(t)$ belonging to the class NS1, if
$S_x(\omega)$ is bandlimited with bandwidth $B$, then $x(t)$  admits
the following mean-square equivalent ``sample representation" \bea
E[( x(t)- \sum_{n=-\infty}^\infty x(nT_s) \frac{sin
\frac{\pi(t-nT_s)}{T_s}}{\frac{\pi(t-nT_s)}{T_s}})^2]  = 0, \
\forall \ t \eea if $T_s < 1/(2B)$,
i.e. $x(t)$ can be reconstructed ``exactly" (in the mean square sense) from its samples taken less than $1/(2B)$ interval apart.
\end{corollary}
\begin{figure}[b!]
\centerline{
\includegraphics
[height=5.5cm,width=9cm]{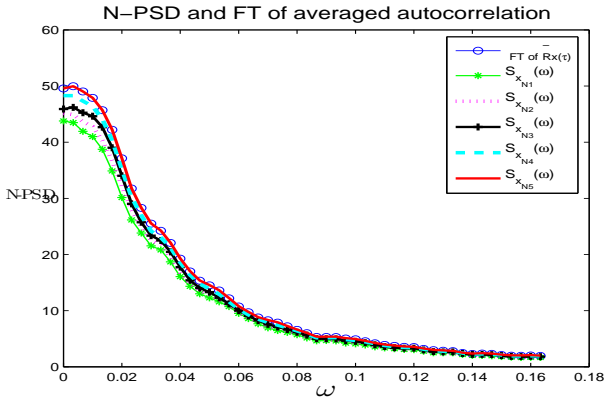}} \vspace{-0.2in}
\caption{{\small Demonstrating Theorem 1 for a nonstationary white noise passed through an IIR filter (explained in the text). 
The blue `-o' line is the Fourier transform of averaged autocorrelation (RHS of
(\ref{nswk})) with the averaged autocorrelation approximated using
$N=500$-length signals. The other five lines are the $N_i$-PSD's
with $N_i$ increasing from $N_1=400$ to $N_5=500$. As can be seen,
this approaches RHS as $N$ increases (here since RHS is approximated
using 500 points, LHS is exactly equal to RHS for $N_5=500$).
   }}
\label{wkverify}
\end{figure}

\section{Application to Time Sequences of Spatially Non-WSS Signals}

To compute the PSD and its bandwidth for non-WSS processes, multiple
realizations of the process are required. One application obtaining these is to deal with a time sequence of
spatial signals. We can compute an estimate of the spatial PSD using
the signals from a temporally i.i.d or a temporally stationary and
ergodic sequence as the multiple realizations. In this section, we
demonstrate Theorem 1 and Corollary 1.

In all simulations, we use discrete time signals. The spatial index
is denoted by $k$ and the time index by $n$, i.e. $x_n[k]$,  as a
function of $k$, is a spatial signal. Averaging over time, $n$, is
used to compute expectations.
For Fig.\ref{wkverify}, 
we generate $x_n[k]$ by passing temporally i.i.d. spatially nonstationary white noise, $w_n[k]$, with $\sigma_{w,n}^2[1:N] = [1_{N/3} , 0.1_{2N/3}]$ (where $c_j$ denotes a vector of $c$'s of length $j$), through a stable infinite impulse response (IIR) filter, $h[k]$, with transfer function $H(\omega) = 1/(1-a_1 e^{-j\omega} - a_2 e^{-2j\omega})$ and $a_1=0.8$, $a_2=0.1$. This is done for each $n$, i.e. $x_n[k] = w_n[k] \star h[k]$, $\forall n$. The signals, $x_n[k]$, belong to the class NS1 and so they satisfy the assumptions of Theorem 1 (as explained in Section 3). We show the verification of Theorem 1 in Fig. \ref{wkverify}. The RHS of (\ref{nswk}) (Fourier transform of averaged autocorrelation function) is plotted as a blue `-o' line. Autocorrelation and its spatial average are computed using an $N=500$ length signal (to approximate $\infty$). We also plot the $N$-length PSD, i.e. the expectation (here average over time) of the normalized periodogram of an $N$-length signal, for increasing $N$ values. 

\begin{figure}[t!] 
\includegraphics
[width=8.5cm, height=9.5cm]{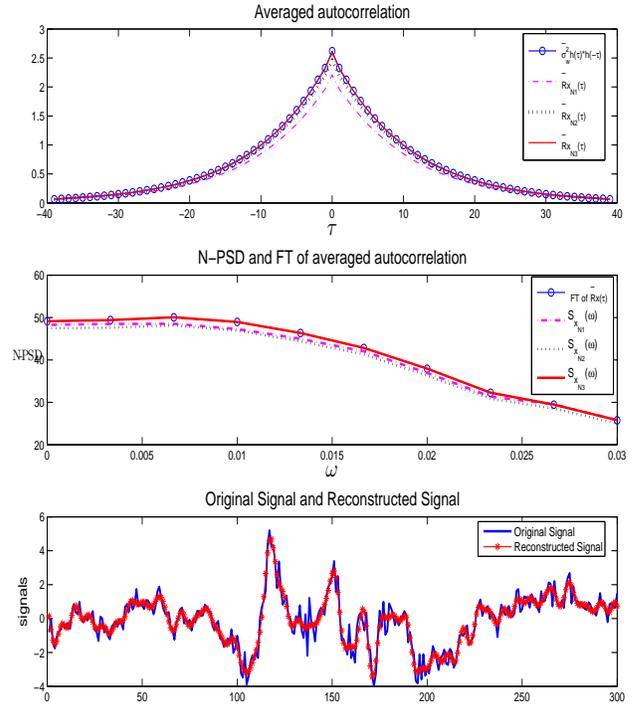} \vspace{-0.2in}
  \caption{\small{Demonstrating parts of Corollary 1.
  These figures are for temporally i.i.d spatially nonstationary white noise passing through a stable IIR AR-2 filter with $H(\omega)=1/ (1-a_1e^{-j\omega}-a_2e^{-j2\omega})$.
The first figure verifies (\ref{rxbar_ns1}). The middle figure
demonstrates (\ref{SxNS1}) and in fact also Theorem 1.  In the
bottom one, one spatial signal and its reconstruction by subsampling
using a 90\%-bandwidth (computed using the N=500-PSD) and
reconstructing the signal is shown.The  mean of the square of the
residual error over all sequences is 8.05\% (will approach 10\% if
the $N$ in the PSD is increased to $\infty$).}}
\label{nyquistverify}
\end{figure}

In Fig. \ref{nyquistverify}, we demonstrate Corollary 1 for the same data, but now using the knowledge of $\sigma_{w,n}^2[k]$ and $h[k]$. The theoretically computed averaged autocorrelation, $\bar{\sigma_w^2}h(\tau) \star h(-\tau)$, and its Fourier transform, $\bar{\sigma_w^2}|H(w)|^2$, are plotted as a blue `-o' line  in the top two figures. The numerically computed averaged autocorrelation and $N$-PSD are also plotted for increasing $N$ in the same two figures. Next, we use the $N$-PSD, for $N=500$, to compute the 90\%-bandwidth (point beyond which the residual PSD sum is less than 10\% of the total PSD sum)  and decimate all $x_n$'s at twice this rate. The computed MSE between $x_n$ and its reconstruction, $\hat{x}_n$, is 8.05\%. This MSE approaches 10\% as $N \tends \infty$. 
Similar results are obtained when we simulated $x_n[k]$ from a temporally stationary process instead of a temporally i.i.d. one. This is done by generating each $w_n[k]$ from a temporal AR-1 process and passing it through the same IIR $h$ as in Fig. \ref{nyquistverify}.

\section{Conclusions}
We have proved the Wiener-Khinchin result for non-WSS processes. This has been combined with Gardner's result \cite{samplenonstat} to prove that Nyquist's criterion can be used to subsample a certain class of PSD-bandlimited nonstationary signals. Application of these two results to subsampling a simulated time sequence of spatially non-WSS signals is shown. Future work includes proving a general Nyquist-type result for all non-WSS signals that satisfy the assumptions of Theorem \ref{wk}.

\end{document}